\documentclass{amsart}

\newtheorem{theorem}{Theorem}
\newcommand{\bt}{\begin{theorem}}
\newcommand{\et}{\end{theorem}}
\newtheorem{lemma}{Lemma}
\newcommand{\bl}{\begin{lemma}}
\newcommand{\el}{\end{lemma}}
\newtheorem{corollary}{Corollary}
\newcommand{\bc}{\begin{corollary}}
\newcommand{\ec}{\end{corollary}}
\newcommand{\beq}{\begin{equation}}
\newcommand{\eeq}{\end{equation}}
\newcommand{\benum}{\begin{enumerate}}
\newcommand{\eenum}{\end{enumerate}}

\newcommand{\Z}{\ensuremath{\mathbf Z}}

\newcommand{\R}{\ensuremath{\mathbf R}}

\usepackage{amsmath,amssymb,graphicx}

\title{Pick's theorem and sums of lattice points}
\author{Karl Levy}
\address{Department of Mathematics, Borough of Manhattan Community College (CUNY), New York, NY 10007}
\email{karl.ethan.levy@gmail.com}
\author{Melvyn B. Nathanson}
\address{Department of Mathematics, Lehman College (CUNY), Bronx, NY 10468}
\email{melvyn.nathanson@lehman.cuny.edu}

\date{}

\begin{document}

\maketitle

\begin{abstract}
Pick's theorem is used to prove that if $P$ is a lattice polygon (that is, the convex 
hull of a finite set of lattice points in the plane), then every lattice point 
in the $h$-fold sumset $hP$ is the sum of $h$ lattice points in $P$.
\end{abstract}

For  sets $X$ and $Y$ in $\R^n$, we define the \emph{sumset} 
\[
X+Y = \{x + y: x \in X \text{ and } y \in Y\}.
\]
For every positive integer $h$, 
we have the \emph{$h$-fold sumset}
\[
hX = \{x_1 + \cdots + x_h: x_i \in X \text{ for } i=1,\ldots, h\}
\]
and the \emph{dilation}  
\[
h \ast X = \{hx:x \in X\}.
\]

\bl
For every convex subset of $\R^n$ and every positive integer $h$, 
the sumset equals the dilation, that is,
\[
hX = h \ast X.
\]
\el

\begin{proof}
If $x \in X$, then 
\[
hx = \underbrace{x + \cdots + x}_{\text{$h$ summands} } \in hX
\]
and so $h\ast X \subseteq hX$.

If  the set $X$ is convex and if $x_1,\ldots, x_h \in X$, then 
$X$ contains the convex combination 
\[
(1/h)x_1 + \cdots + (1/h)x_h 
\]
and so 
\[
x_1 + \cdots + x_h = h \left( (1/h)x_1 + \cdots + (1/h)x_h  \right) \in h\ast X.
\]
Thus, if $X$ is convex, then $hX \subseteq h\ast X$.  
This completes the proof. 
\end{proof}

A \emph{lattice polytope}  in $\R^n$ is the convex hull of a nonempty 
finite set of lattice points  in $\Z^n$.  
A \emph{lattice polygon}  in $\R^2$ is the convex hull of a nonempty finite subset of $\Z^2$.  

Let $n \geq 2$ and $h \geq 2$.
If $P, P_1, \ldots, P_h$ are lattice polytopes in $\R^n$, then
\beq            \label{PolySum:hP}
h(P \cap \Z^n) \subseteq (hP) \cap \Z^n
\eeq
and
\beq            \label{PolySum:P1Ph}
(P_1 \cap \Z^n) + \cdots + (P_h\cap \Z^n) \subseteq (P_1 + \cdots + P_h) \cap \Z^n.  
\eeq
These set inclusions can be strict.  For example, 
if $P_1$ is the triangle in $\R^2$ whose vertices are $\{(0,0), (1,0), (1,-1)\}$, 
and if $P_2$ is the triangle whose vertices are $\{(0,0), (1,2), (2,3)\}$, 
then $P_1+P_2$ is the hexagon with vertices 
\[
\{(0,0),  (1,-1), (3,2), (3,3), (2,3), (1,2) \}.
\]
We have
\[
(1,1) = (1/2,0) + (1/2,1)  \in (P_1+P_2)\cap \Z^2
\]
but 
\[
(1,1) \notin (P_1\cap \Z^2)  + (P_2\cap \Z^2).
\]
Therefore,
\[
(P_1 \cap \Z^2) + (P_2\cap \Z^2) \neq (P_1 + P_2) \cap \Z^2.  
\]

In $\R^3$, let $P$ be the tetrahedron with vertices 
\[
\{ (0,0,0), (1,0,0), (0,1,0), (1,1,2) \}
\]
We have
\[
(1,1,1) = (1/2,1/2,0) + (1/2, 1/2,1) \in (2P)\cap \Z^3
\]
but 
\[
(1,1,1)\notin 2(P \cap \Z^3).
\]
Therefore,
\[
2(P \cap \Z^3)  \neq (2P) \cap \Z^3.  
\]

It is an only partially solved problem to determine the lattice polytopes 
$P$, $P_1$, \ldots, $P_h$ in $\R^n$ for which 
we have equalities and not inclusions in~\eqref{PolySum:hP}
and~\eqref{PolySum:P1Ph}.  
This is usually discussed in the language of toric geometry~\cite{fakh02,haas08,oda08}.
It is known that in the plane we have 
\beq   \label{PolySum:BIG}
h(P \cap \Z^2) = (hP) \cap \Z^2
\eeq
for every lattice polygon $P$ and every positive integer $h$
(Koelman~\cite{koel93a,koel93b}).  
In this note we apply Pick's theorem (Pick~\cite{pick99}, 
Beck and Robins~\cite[pp.  38--40]{beck-robi07}) to obtain a simple proof of this result.  

 Pick's theorem, proofs of which appear frequently in the \emph{Monthly} 
(e.g.~\cite{funk74,varb85,diaz-robi95,murt-thai07}), 
states that if $P$ is a lattice polygon with area $A$ 
and with $I$ lattice points in its interior and $B$ lattice points 
on its boundary, then 
\beq   \label{PolySum:Pick}
A = I + \frac{B}{2}  - 1.
\eeq
A triangle is the convex hull of three non-collinear points.  
A \emph{primitive lattice triangle} is a lattice triangle whose only lattice points 
are its three vertices.  
By Pick's theorem, a lattice triangle in $\R^2$ is primitive if and only if its area is $1/2$.

\bl               \label{PolySum:lemma:PrimitiveTriangle}
Let $T$ be a primitive lattice triangle with vertex set $\{0, u,v\}$, and let  
\[
W = \{ iu + jv : i = 0,1,\ldots, h \text{ and } j = 0,1,2,\ldots, h-i\}.
\] 
Then
\[
W = hT \cap \Z^2.
\]
\el

\begin{proof}
Because $u,v \in \Z^2$ and because  the set  $T$ is convex and 
\[
hT = h\ast T = \{ \alpha u + \beta v : \alpha \geq 0, \beta \geq 0, \text{ and } \alpha + \beta \leq h \}
\]
it follows that  $W \subseteq hT\cap \Z^2$.  
The linear independence of the vectors $u$ and $v$ implies that   
\[
|W| = \sum_{i=0}^h (h-i+1) = \frac{(h+1)(h+2)}{2}.
\]

Let $A(h)$ denote the area of the dilated triangle $hT$, and let $I(h)$ and $B(h)$ denote, 
respectively,  the number of interior lattice points and boundary lattice points of $hT$. 
Because $T$ is primitive, we have $A(1) = 1/2$ and $B(1) = 3$.  
It follows that  
\[
A(h) =  A(1) h^2= \frac{h^2}{2}
\]
and
\[
B(h) =  B(1)h = 3h.
\]
By Pick's theorem, 
\[
A(h) = I(h) + \frac{B(h)}{2}-1
\]
and so the number of lattice points in $hT$ is 
\[
|hT \cap \Z^2 | = I(h) + B(h) = A(h) + \frac{B(h)}{2}+1 = \frac{(h+1)(h+2)}{2}.
\]
Because $W$ and $hT\cap \Z^2$ are finite sets with 
$W \subseteq hT\cap \Z^2$ and $|W| = |hT\cap \Z^2|$, it follows that 
$W = hT\cap \Z^2$.  This completes the proof.  
\end{proof}

\bt
Let $P$ be a lattice polygon.  
If $w$ is a lattice point in the sumset $hP$, then there exist 
lattice points $a,b,c$ in $P$ and nonnegative integers $i$, $j$, and $k$ 
such that $h = i+j + k$ and 
\[
w = ia + jb + kc \in  h(P\cap \Z^2).
\]
In particular, 
\[
h(P\cap \Z^2) = (hP)\cap \Z^2.  
\]
\et

\begin{proof}
Every lattice polygon $P$ can be triangulated into primitive lattice triangles.  
If $w$ is a lattice point in $hP = h\ast P$, then $w/h \in P$ and so there is a primitive lattice 
triangle $T'$ contained in $P$ with $w/h \in T'$.  Let $\{a,b,c\}$ be the set of vertices of $T'$, 
and let $T = T'-c$ be the primitive lattice triangle with vertices $0$, $u = a-c$, and $v = b-c$.
We have $w/h - c \in T$ and so $w-hc \in hT$.  
By Lemma~\ref{PolySum:lemma:PrimitiveTriangle}, there are nonnegative integers 
$i$ and $j$ such that $i+j = h-k \leq h$ and 
\[
w - hc = iu + jv.
\]
This implies that  
\[
w = iu + jv + hc = ia + jb + kc \in  3(T' \cap \Z^2) \subseteq 3(P\cap \Z^2).
\]
This completes the proof.  
\end{proof}

\def\cprime{$'$} \def\cprime{$'$} \def\cprime{$'$}
\providecommand{\bysame}{\leavevmode\hbox to3em{\hrulefill}\thinspace}
\providecommand{\MR}{\relax\ifhmode\unskip\space\fi MR }
\providecommand{\MRhref}[2]{%
  \href{http://www.ams.org/mathscinet-getitem?mr=#1}{#2}
}
\providecommand{\href}[2]{#2}


\begin{thebibliography}{10}

\bibitem{beck-robi07}
M.~Beck and S.~Robins, \emph{{Computing the Continuous Discretely}},
  Springer, New York, 2007.

\bibitem{diaz-robi95}
R.~Diaz and S.~Robins, \emph{Pick's formula via the {W}eierstrass
  {$\wp$}-function}, Amer. Math. Monthly \textbf{102} (1995), no.~5, 431--437.

\bibitem{fakh02}
N.~Fakhruddin, \emph{{Multiplication maps of linear systems on smooth
  projective toric varieties}}, {arXiv: 0208178}, 2002.

\bibitem{funk74}
W.~W. Funkenbusch, \emph{Classroom {N}otes: {F}rom {E}uler's {F}ormula to
  {P}ick's {F}ormula {U}sing an {E}dge {T}heorem}, Amer. Math. Monthly
  \textbf{81} (1974), no.~6, 647--648.

\bibitem{haas08}
C.~Haase, B.~Nill, A.~Paffenholz, and F.~Santos, \emph{Lattice points in
  {M}inkowski sums}, Electron. J. Combin. \textbf{15} (2008), no.~1, Note 11,
  5.

\bibitem{koel93b}
R.~J. Koelman, \emph{A criterion for the ideal of a projectively embedded toric
  surface to be generated by quadrics}, Beitr\"age Algebra Geom. \textbf{34}
  (1993), no.~1, 57--62.

\bibitem{koel93a}
\bysame, \emph{Generators for the ideal of a projectively embedded toric
  surface}, Tohoku Math. J. (2) \textbf{45} (1993), no.~3, 385--392.

\bibitem{murt-thai07}
M.~Ram Murty and N.~Thain, \emph{Pick's theorem via {M}inkowski's theorem},
  Amer. Math. Monthly \textbf{114} (2007), no.~8, 732--736.

\bibitem{oda08}
T.~Oda, \emph{{Problems on Minkowski sums of convex lattice polytopes}},
  {arXiv: 0812.1418}, 2008.

\bibitem{pick99}
G.~A. Pick, \emph{{Geometrisches zur Zahlenlehre}}, {Sitzenber. Lotos (Prague)}
  \textbf{19} (1899), 311--319.

\bibitem{varb85}
D.E. Varberg, \emph{Pick's theorem revisited}, Amer. Math. Monthly \textbf{92}
  (1985), no.~8, 584--587.

\end{thebibliography}
\end{document}